\documentclass[11pt]{amsart}
\usepackage{amsmath,amsfonts,latexsym,amssymb,amscd, graphicx,pictexwd}

\renewcommand{\P}{{\mathbb P}}
\newcommand{\E}{{\mathbb E}}
\newcommand{\R}{{\mathbb R}}

\newcommand{\ZZ}{{\mathbb Z}}
\newcommand{\PP}{{\mathbf P}}

\newcommand{\0}{{\mathbf 0}}
\newcommand{\pp}{{\mathbf p}}
\newcommand{\qq}{{\mathbf q}}
\newcommand{\cc}{{\mathbf c}}
\newcommand{\zz}{{\mathbf z}}
\newcommand{\yy}{{\mathbf y}}
\newcommand{\mm}{{\mathbf m}}
\newcommand{\ww}{{\mathbf w}}
\newcommand{\xx}{{\mathbf x}}
\newcommand{\nn}{{\mathbf n}}
\newcommand{\cal}{\mathcal}

\newcommand{\calN}{{\mathcal N}}

\newcommand{\Z}{{\mathbb Z}}
\newcommand{\calF}{{\mathcal F}}

\newcommand{\eps}{{\varepsilon }}

\newcommand{\argmax}{\mathop{\rm argmax}}

\newtheorem{thm}{Theorem}[section]
\newtheorem{coro}[thm]{Corollary}
\newtheorem{lem}[thm]{Lemma}
\newtheorem{prop}[thm]{Proposition}

\numberwithin{equation}{section}

\begin{document}
\title[Busemann functions and Equilibrium measure]{Busemann functions and Equilibrium measures\\ in last passage percolation models}

\author{Eric Cator}
\address{Delft University of Technology\\
Mekelweg 4, 2628 CD Delft, The Netherlands}
\email{E.A.Cator@tudelft.nl}

\author{Leandro P. R. Pimentel}
\address{Institute of Mathematics\\
Federal University of Rio de Janeiro}
\email{leandro@im.ufrj.br}
\thanks{Leandro P. R. Pimentel was supported by grant number 613.000.605 from the Netherlands Organisation for Scientific Research (NWO)}

\begin{abstract}
The interplay between two-dimensional percolation growth models and one-dimensional particle processes has been a fruitful source of interesting mathematical phenomena. In this paper we develop a connection between the construction of Busemann functions in the Hammersley last-passage percolation model with i.i.d. random weights, and the existence, ergodicity and uniqueness of equilibrium (or time-invariant) measures  for the related (multi-class) interacting fluid system. As we shall see, in the classical Hammersley model, where each point has weight one, this approach brings a new and rather geometrical solution of the longest increasing subsequence problem, as well as a central limit theorem for the Busemann function.
\end{abstract}

\maketitle

\section{Introduction}\label{sec:intro}

In the middle of the fifties H. Busemann \cite{Bu} introduced a collection of functions to study geometrical aspects  of metric spaces. These functions are induced by a metric $d$, and by a collection of rays (semi-infinite geodesics) as follows: the Busemann function $b_\varpi(\cdot)$, with respect to a ray $(\varpi(r))_{r\geq 0}$, is the limit of $d(\varpi(r),\varpi(0))-d(\varpi(r),\cdot)$ as $r$ goes to infinity. Along a ray $\varpi$ the metric $d$ becomes additive. By using the triangle inequality, this implies that the defining sequence is nondecreasing and bounded from above, and so it always converges. Using analogous considerations, one can construct Busemann functions over spaces equipped with a super-additve ``metric''  $L$ (one needs the reversed triangle inequality). In this work we are particularly interested in geometrical aspects of the following stochastic  two-dimensional last passage (super-additive) percolation model: let $\PP\subseteq\R^2$ be a two-dimensional Poisson random set of intensity one. On each point $\pp\in\PP$ we put a random positive weight $\omega_\pp$ and we assume that $\{\omega_\pp\,:\,\pp\in\PP\}$ is a collection of i.i.d. random variables, distributed according to a distribution function $F$, which are also independent of $\PP$. When $F$ is the Dirac distribution concentrated on $1$ (each point has weight $1$; we will denote this $F$ by $\delta_1$), then we refer to this model as the classical Hammersley model (Aldous \& Diaconis \cite{AlDi}). For each $\pp,\qq\in\R^2$, with $\pp<\qq$ (inequality in each coordinate, $\pp\neq \qq$), let $\Pi(\pp,\qq)$ denote the set of all increasing (or up-right) paths, consisting of points in $\PP$, from $\pp$ to $\qq$, where we exclude the starting point $\pp$. In this probabilistic model, the ``metric'' (or last-passage time) $L$ between  $\pp\leq\qq$ is defined by
$$L(\pp,\qq):=\max_{\varpi\in\Pi(\pp,\qq)}\big\{\sum_{\pp'\in\varpi}\omega_{\pp'}\big\}\,.$$
Then $L$ is super-additive,
$$L(\pp,\qq)\geq L(\pp,\zz)+L(\zz,\qq)\,.$$
When we consider a path $\varpi$ from $\pp$ to $\qq$ consisting of increasing points $(\pp_1,\ldots,\pp_n)$, we will view $\varpi$ as the lowest increasing continuous path connecting all the points, starting at $\pp$ and ending at $\qq$, and then excluding $\pp$. This way we can talk about crossings with other paths or with lines. A finite geodesic between $\pp$ and $\qq$ is given by the lowest path that attains the maximum in the definition of $L(\pp,\qq)$, which we will denote by $\varpi(\pp,\qq)$ (this is well defined for any ordered pair $(\pp,\qq)$, even if we do not specify the order).\\

In \cite{CPshape}, using methods developed by Newman and co-authors in \cite{Ne} and \cite{HoNe} and also applied to the classical Hammersley process by W\"uthrich in \cite{Wu}, it is shown that these finite geodesics can be extended to semi-infinite $\alpha${\em -rays} by moving one endpoint appropriately to infinity. An $\alpha$-ray starting at $\pp\in\R^2$, denoted by $\varpi_\alpha(\pp)$, is a semi-infinite geodesic that starts at $\pp$ and moves to infinity in the direction $\vec{\alpha}:=(\cos\alpha,\sin\alpha)$, for $\alpha\in (\pi, 3\pi/2)$. It turns out that for fixed $\alpha$, with probability $1$, each $\pp\in\R^2$ is the starting point of a unique $\alpha$-ray, and two $\alpha$-rays will always coalesce eventually. In Section \ref{sec:alpharay} we will state precisely the theorems for existence and coalescence of $\alpha$-rays. An important tool in the development of these concepts is the following result, known as the shape theorem: there exists a constant $\gamma=\gamma(F)$ such that for all $x,t\geq 0$
\begin{equation}\label{eq:introshape}
\lim_{r\to \infty} \frac{L(\0,r(x,t))}{r} = \gamma\sqrt{xt}\ \ ({\rm a.s.})
\end{equation}

In Section \ref{sec:busemann} we will use the $\alpha$-rays to construct the Busemann function $B_\alpha$: if $\xx, \yy \in \R^2$, let $\cc=\cc(\alpha,\xx,\yy)$ be the coalescence point of the two $\alpha$-rays starting at $\xx$ and $\yy$. Then
\[ B_\alpha(\xx,\yy) = L(\cc,\yy) - L(\cc,\xx).\]
It is not hard to see that the distribution of $B_\alpha$ is invariant under translation and that $B_\alpha$ is additive: $B_\alpha(\xx,\zz) = B_\alpha(\xx,\yy) + B_\alpha(\yy,\zz)$. Its most important property, however, is the following connection with an associated interacting particle process: if we define the measures $\nu_\alpha^t$ on $\R$ by
\[ \nu_\alpha^t((x,y]) = B_\alpha((x,t),(y,t)),\]
then the family of random measures $\{\nu_\alpha^t\ :\ t\in \R\}$ forms a Markov process, and its evolution corresponds to the evolution of an interacting particle process (maybe interacting {\em fluid} process is a better name) which is a natural extension of the classical Hammersley interacting particle process, where $F=\delta_1$. This connection turns out to be the key idea of this paper.  On the one hand, it allows us to prove existence, uniqueness and ergodicity (mixing) of the equilibrium measures of the interacting fluid process. Furthermore, it gives us a natural way to prove a strong law of large numbers for the second class particle and to define a multi-class fluid system with a countable number of classes, which is a new result even in the classical Hammersley process. Also, the classical result that $\gamma(\delta_1)=2$ (\eqref{eq:introshape} for the classical Hammersley process) follows easily from our methods. On the other hand, it implies a central limit theorem for the Busemann function in the classical model. This result also shows a phase transition from a Gaussian limit distribution, on the square-root scale, to (zero-mean)  Tracy-Widom type limit distribution, on the cube-root scale, at the critical angle $\alpha-\pi$. This transition from square-root to cube-root scaling for Busemann functions was first conjectured in 2001 by Howard \& Newman \cite{HoNe} (see also W\"utrich \cite{Wu}).\\

In Section \ref{sec:TASEP} we describe how to define a Busemann function in the last passage percolation on the lattice $\Z^2$ with iid weights on the lattice points. It is an important open question for a long time already, how to prove in this general setup that the shape function, see \eqref{eq:introshape}, is strictly curved, a fact we need to define our $\alpha$-rays. However, if we restrict ourselves to exponential (or geometric) weights, in which case the last passage percolation is an alternative description of the totally asymetric exclusion process (TASEP), we know the shape function, and it is indeed strictly curved. Therefore, we can define the Busemann function and we find a similar connection to the (known) equilibrium measures, which allows us to prove analogous results.\\

We feel that the connection between the Busemann functions and the equilibrium measures gives us an important new tool to study last passage percolation and the corresponding fluid processes. In an upcoming paper, we will use the Busemann function and the results from this paper to determine the asymptotic speed of a second class particle, given a deterministic, rarefaction initial condition, in the classical Hammersley process and in TASEP. Furthermore, we have strong indications that the Busemann function can help us establish the cube-root behavior of the length of a longest path, and the fluctuations of the longest path.\\

\paragraph{\bf Overview} In Section \ref{sec:alpharay} we will state the theorems, that we will use further on, about the limit shape and the existence and coalescence of $\alpha$-rays, and then give precise references for the proofs. In Section \ref{sec:busemann} we define the Busemann function and give its most important properties. In Section \ref{sec:fluidprocess} we introduce the Hammersley interacting fluid process and establish the connection between the Busemann function and the equilibrium measures. In Section \ref{sec:ergodic} we prove uniqueness and ergodicity (mixing property) of the equilibrium measures, and we show local convergence to the equilibrium measure in case of a rarefaction fan. In Section \ref{sec:scalBus} we prove the central limit theorem for the Busemann function in the classical model. In Section \ref{sec:multiclass}, we show how we can define a multi-class system with a countable number of classes, and we establish the strong law for a second class particle. In Section \ref{sec:TASEP} we state the analogous results for the TASEP.

\section{Shape function and $\alpha$-rays}\label{sec:alpharay}
A key notion in this paper will be an \emph{$\alpha$-ray}: for each angle $\alpha\in (\pi, 3\pi/2)$ and for each point $\xx\in\R^2$, $\varpi_\alpha(\xx)$ is the lowest continuous down-left path through an ordered sequence $(\pp_i)_{i\geq 0}$ in $\R^2$, with $\pp_0=\xx$, $\pp_i\in \PP$ and $\pp_i\geq \pp_j$ whenever $ i\leq j$. Furthermore, $\varpi(\pp_j,\pp_i) \subset \varpi_\alpha(\xx)$  (every  part of the path is a geodesic), and finally we must have that
\begin{equation}\label{eq:ray}
 \lim_{i\to \infty} \frac{\pp_i}{\|\pp_i\|} = \vec{\alpha}:=(\cos\alpha, \sin\alpha)\,.
\end{equation}
A crucial step for the existence of $\alpha$-rays is the following shape theorem: set $\0=(0,0)$, $\nn=(n,n)$,
\[F(x) = \P(\omega_\pp \leq x)\,\mbox{ and }\,\gamma=\gamma(F)=\sup_{n\geq 1}\frac{\E(L(\0,\nn))}{n}> 0 \,.\]
\begin{thm}\label{thm:shape}
Suppose that
\begin{equation}\label{eq:a1}
\int_0^\infty \sqrt{1-F(x)}\,dx <+\infty\,.
\end{equation}
Then $\gamma(F)<\infty$ and for all $x,t>0$, as $r \to \infty$,
\[\frac{L\left(\0,(rx,rt)\right)}{r}\to \gamma \sqrt{xt}\ \  {\rm a.s.}\ \ \ \mbox{and} \ \ \   \frac{\E L\left(\0,(rx,rt)\right)}{r} \to \gamma \sqrt{xt}\,.\]
Further, if \eqref{eq:a1} is strengthened to: there exists $a>0$ such that
\begin{equation}\label{eq:a2}
 \int_0^\infty \exp(a x)\,dF(x) <+\infty\,,
\end{equation}
then there exist constants $c_0,c_1,c_2,c_3,c_4>0$ such that for all $r\geq c_0$
\[\P\big(|L(\0,(r,r))-\gamma r |\geq u\big)\leq c_1\exp\Big(-c_2\frac{u}{\sqrt{r}\log r}\Big)\,\]
for  $u\in \left[\,c_3 \sqrt{r}\log^2 r\,,\,c_4 r^{3/2}\log r\,\right]$.
\end{thm}
For proofs see Theorem 1.1 and 1.2 in \cite{CPshape}. Theorem \ref{thm:shape} shows that $L$ has a curved limiting shape, mainly due to the invariance of the Poisson process under volume preserving maps: if $x,t,r >0$ and $\pp\in \R^2$, then
\begin{equation}\label{eq:sym}
L(\0,(x,t)) \stackrel{\cal D}{=} L\left(\pp,\pp + (r x,t/r)\right)\,.
\end{equation}
This is because under this map, the distribution of the Poisson process does not change, and the up-right paths are preserved. The almost sure convergence is a standard consequence of the sub-additive ergodic theorem, once we have a bound on $\E(L(\0,\nn))$, linear in $n$.\\

The following theorem gives us existence and coalescence of $\alpha$-rays. Before we state the theorem, we shall define what we mean by convergence of paths: we say that a sequence of paths $\varpi^n$ converges to $\varpi$, and denote
$\lim_{n\to\infty}\varpi^n= \varpi$, if for all bounded subsets  $B\subset \R^2$ there exists $n_0$ such that  $\varpi^n\cap B=\varpi\cap B$ for all $n\geq n_0$.
\begin{thm}\label{thm:existray}
Assume \eqref{eq:a2}, so for some $a>0$,
\begin{equation}
\nonumber \int_0^\infty \exp(a x)\,dF(x) <+\infty\,.
\end{equation}
Then for fixed $\alpha\in(\pi,3\pi/2)$ the following holds with probability one:
\begin{enumerate}
\item For each $\xx\in\R^2$ there exists a unique $\alpha$-ray starting from $\xx$, which we denote by $\varpi_\alpha(\xx)$.
\item For any sequence $(\zz_n)_{n\geq 0}$ of points in $\R^2$ with $\|\zz_n\|\to \infty$,
$$\mbox{if }\lim_{n\to \infty} \frac{\zz_n}{\|\zz_n\|} = (\cos\alpha, \sin\alpha)\,\mbox{ then }\,\lim_{n\to\infty}\varpi(\xx,\zz_n)=\varpi_\alpha(\xx)\,.$$
\item For all $\xx,\yy\in\R^2$ there exists $\cc=\cc(\alpha,\xx,\yy)$ such that $\varpi_\alpha(\xx)$ and $\varpi_\alpha(\yy)$ coalesce at $\cc$: $$\varpi_\alpha(\xx)=\varpi(\xx,\cc)\cup\varpi_\alpha(\cc)\mbox{ and }\varpi_\alpha(\yy)=\varpi(\yy,\cc)\cup\varpi_\alpha(\cc)\,.$$
\end{enumerate}
\end{thm}
The proof of Theorem \ref{thm:existray} is based on a method introduced by Newman \cite{Ne} that can be applied in a wide percolation context. An outline of the proof is given in \cite{CPshape}, Theorem 2.4 and 2.5. A detailed proof of the same theorem, but now restricted to the Hammersley classical model, can be found in \cite{Wu}.

\section{The Busemann function}\label{sec:busemann}

Using the concept of $\alpha$-rays, we will study the function $B_\alpha(\xx,\yy)$, which is defined by taking the first coalescence point $\cc=\cc(\alpha,\xx,\yy)$ between the $\alpha$-ray that starts from $\xx$ and the one that starts from $\yy$ (remember that these two rays coalesce), and setting
\begin{equation}\label{eq:defLa}
B_{\alpha}(\xx,\yy)=L(\cc,\yy)-L(\cc,\xx)\,.
\end{equation}
Note that if we take a different coalescence point $\cc'$, then $\cc\geq\cc'$ and they both lie on a geodesic. Since $L$ is additive on a geodesic, we get $L(\cc',\xx) = L(\cc',\cc) + L(\cc,\xx)$, which shows that the definition of $B_\alpha(\xx,\yy)$ does not depend on the choice of the coalescence point.
Let $(\zz_n)_{n\geq 1}$ be any unbounded decreasing sequence that follows direction $(\cos\alpha,\sin\alpha)$, and let $\cc(\zz_n,\xx,\yy)$ denote the most up-right coalescence point between $\varpi\left(\zz_n,\xx\right)$ and $\varpi\left(\zz_n,\yy\right)$. By Theorem \ref{thm:existray}, with probability one, there exists $n_0>0$ such that
\begin{equation}\label{eq:coal1}
\forall\,n\geq n_0\,\,\,\,\cc(\zz_n,\xx,\yy)=\cc(\alpha,\xx,\yy)\,\,\mbox{ and }\,\,L(\zz_n,\yy) - L(\zz_n,\xx)=B_\alpha(\xx,\yy)\,.
\end{equation}
Therefore, in geometrical terms, $B_\alpha(\xx,\cdot)$ can be seen as the \emph{Busemann function} along the ray $\varpi_\alpha(\xx)$.\\

Some properties of $B_\alpha$ are summarized in the following proposition. The proofs are relatively straightforward, and can be found in the Appendix.
\begin{prop}\label{prop:Busemann}
Define the Busemann function $B_\alpha$ as above, for $\alpha\in (\pi,3\pi/2)$.
\begin{enumerate}
\item\label{prop:Buse1} The distribution of the function $B_\alpha$ is translation invariant: $\forall\ \pp\in\R^2$ \[B_{\alpha}(\cdot+\pp,\cdot+\pp)\stackrel{\cal D}{=}B_{\alpha}(\cdot,\cdot).\]
\item\label{prop:Buse2} $B_\alpha$ is anti-symmetric and additive: $\forall\ \xx, \yy, \zz \in \R^2$
     \[B_{\alpha}(\xx,\yy)=-B_{\alpha}(\yy,\xx)\ \ \mbox{and}\ \ \ B_\alpha(\xx,\zz) = B_\alpha(\xx,\yy) + B_\alpha(\yy,\zz).\]
\item\label{prop:Busesym} For any $(x,t)\in \R^2$,
\[ B_{5\pi/2-\alpha}(\0,(x,t)) \stackrel{\cal D}{=} B_\alpha(\0,(t,x)).\]
\item\label{prop:Buse3} If $\xx\leq \yy$ and $\xx\neq\yy$, then
\[ B_\alpha(\xx,\yy)\geq 0\ \ \mbox{and}\ \ \ 0<\E(B_\alpha(\xx,\yy))<+\infty.\]
\item\label{prop:Buse4} Fix $\xx, \yy\in\R^2$ and $\pp,\qq\in \R^2$ such that $\pp,\qq\geq\0$. The function $\lambda\mapsto B_\alpha(\xx+\lambda\pp,\yy+\lambda\qq)$ is c\`adl\`ag in $\lambda\in \R$.
\end{enumerate}
\end{prop}

As mentioned in the Introduction, the most important aspect of $B_\alpha$ is a Markovian structure described in the following proposition.
\begin{prop}\label{prop:BusMarkov}
For all $s\leq t$ and $x\in \R$ we have
\[ B_\alpha((0,s),(x,t)) = \sup_{z\leq x}\left\{ B_\alpha((0,s),(z,s)) + L((z,s),(x,t))\right\}.\]
\end{prop}
\noindent{\bf Proof:} Without loss of generality we can take $s=0$ (and therefore $t\geq 0$). Define $Z_\alpha=Z_\alpha(x,t)\in \R$ as the crossing-point of the $\alpha$-ray starting at $(x,t)$ with the $x$-axis. Clearly, $Z_\alpha\leq x$ and
\begin{eqnarray*}
B_\alpha(\0,(x,t)) & = & B_\alpha(\0,(Z_\alpha,0)) + B_\alpha((Z_\alpha,0),(x,t))\\
& = &  B_\alpha(\0,(Z_\alpha,0)) + L((Z_\alpha,0),(x,t)).
\end{eqnarray*}
The last equality follows from the fact that $(x,t)$ and $(Z_\alpha,0)$ are lying on an $\alpha$-ray. This means that it is enough to prove that for all $z\leq x$,
\begin{equation}\label{eq:zZalpha}
B_\alpha(\0,(z,0)) + L((z,0),(x,t))\leq B_\alpha(\0,(Z_\alpha,0)) + L((Z_\alpha,0),(x,t)).
\end{equation}
Suppose $\pp$ is a coalescence point of the $\alpha$-rays starting at $\0$, $(x,t)$ and $(z,0)$. Then
\[ B_\alpha(\0,(z,0)) = L(\pp,(z,0)) - L(\pp,\0)\ \ \ \mbox{and}\ \ \ B_\alpha(\0,(Z_\alpha,0))=L(\pp,(Z_\alpha,0)) - L(\pp,\0).\]
Furthermore, since $\pp, (Z_\alpha,0)$ and $(x,t)$ are elements of $\pi_\alpha(x,t)$, we know that
\begin{eqnarray*}
L(\pp,(z,0)) + L((z,0),(x,t)) & \leq & L(\pp, (x,t))\\
& = & L(\pp,(Z_\alpha,0)) + L((Z_\alpha,0),(x,t)).
\end{eqnarray*}
From this, \eqref{eq:zZalpha} easily follows.

\hfill$\Box$\\

Define the positive measures $\nu_t^\alpha$ on $\R$, using Proposition \ref{prop:Busemann}, by
\begin{equation}\label{eq:defnua}
\nu_t^\alpha((x,y]) := B_\alpha((x,t),(0,0))-B_\alpha((y,t),(0,0))=B_\alpha((x,t),(y,t))\ \ \ \ \ \forall\ x\leq y.
\end{equation}
Proposition \ref{prop:BusMarkov} then shows that the process $t\mapsto \nu_t^\alpha$ is a Markov process: the future evolution of $\nu_t^\alpha$ depends on the Poisson process in the upper-half plane $\R\times (t,\infty)$ and on the present value of $\nu_t^\alpha$, not on the past of the process (which is of course independent of the Poisson process in $\R\times (t,\infty)$). Note that the distribution of $\nu_t^\alpha$ does not depend on $t$ (by Proposition \ref{prop:Busemann}\eqref{prop:Buse1}), so this distribution is an equilibrium (or time invariant) measure for the underlying  Markov process. In the next section we will describe the generator of this Markov process, which will be an extension of the classical Hammersley interacting particle process.

\section{The Hammersley interacting fluid system}\label{sec:fluidprocess}

It is well known that the classical Hammersley model, where all weights are 1, described in Aldous \& Diaconis \cite{AlDi}, has a representation as an interacting particle system. The Hammersley process with random weights has a similar description, although a better name might be an interacting fluid system. We start by restricting the compound Poisson process $\{\omega_\pp\,:\,\pp\in\PP\}$ to $\R\times\R_+$. Then we choose a positive, locally finite measure $\nu$ defined on $\R$. Usually, these measures will be purely atomic, but this is not necessary. To each measure $\nu$ we associate a non-decreasing process $\nu(\cdot)$ defined by
\[\nu(x)=\left\{\begin{array}{ll} \nu([0,x]) & \mbox{for } x\geq 0\\
-\nu((x,0)) & \mbox{for } x<0.\end{array}\right.\]
Note that $\nu(\cdot)$ is a cadlag function. Although the details are a bit cumbersome, all the results we will show can be extended, mutatis mutandis, to the case where $\nu(x)=-\infty$ for $x<0$, which would correspond to a non-locally finite measure with an infinite fluid density to the left of $0$. This is a quite natural starting condition, but we will not use it explicitly in this paper.

\begin{figure}[!ht]
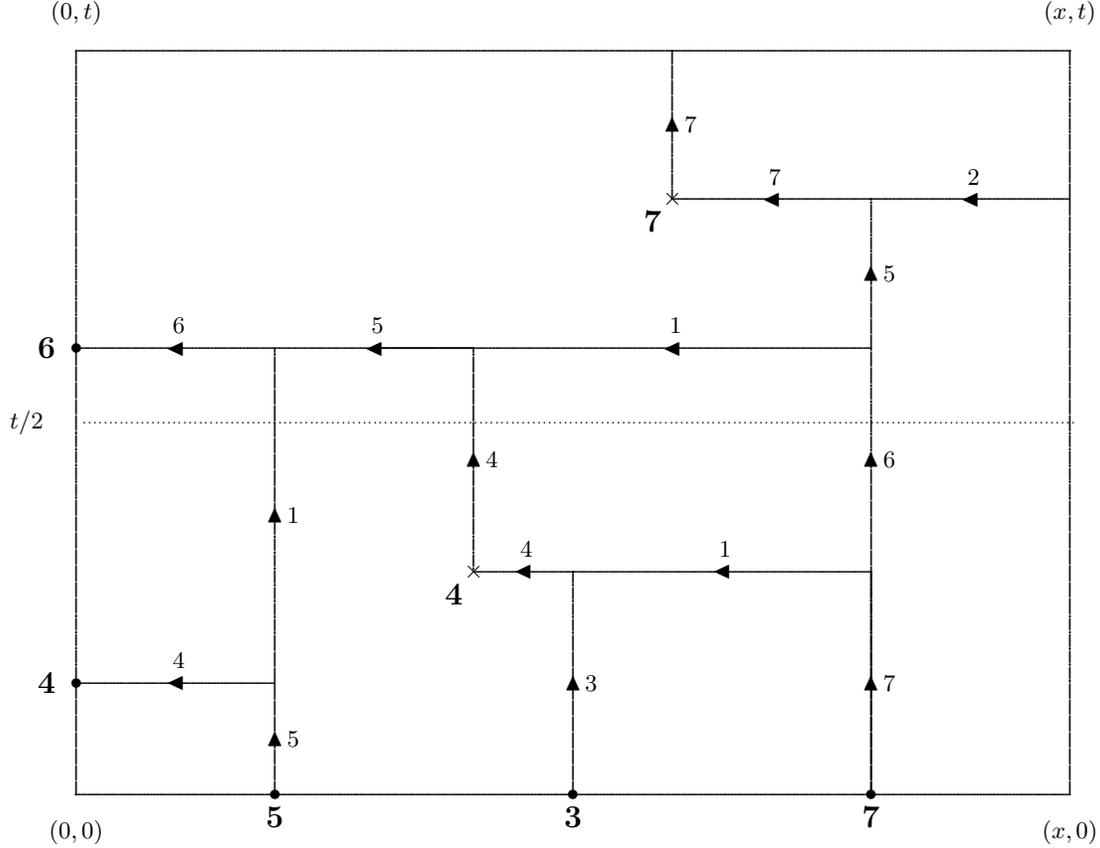

\begin{center}
\strut
$$
\beginpicture
\footnotesize
\setcoordinatesystem units <0.08\textwidth,0.06\textwidth>
  \setplotarea x from 0 to 10, y from 0 to 12
\multiput {\bf $\times$} at
4 3
6 8
/

\multiput {$\bullet$} at
2 0
5 0
8 0
0 1.5
0 6
/
\put {$(0,0)$} at 0.0 -0.5
\put {$(x,0)$} at 10 -0.5
\put {$(0,t)$} at 0.0 10.5
\put {$(x,t)$} at 10 10.5
\put {$t/2$} at -0.5 5

\setlinear
\plot 0 10  10 10  10 0 0 0 0 10
/
\plot 10 8 8 8
/
\plot 8 0 8 8 6 8 6 10
/
\plot 8 0 8 3 5 3
/
\plot 5 0 5 3 4 3 4 6 0 6
/
\plot 8 6 3 6
/
\plot 2 0 2 6
/
\plot 2 1.5 0 1.5
/
\multiput {$\blacktriangle$} at
2 0.75
5 1.5
8 1.5
8 4.5
8 7
6 9
4 4.5
2 3.75
/
\multiput {$\blacktriangleleft$} at
1 1.5
1 6
3 6
6 6
6.5 3
7 8
9 8
4.5 3
/

{\large
\put {\bf 5} at 2 -0.3
\put {\bf 3} at 5 -0.3
\put {\bf 7} at 8 -0.3
\put {\bf 4} at -0.3 1.5
\put {\bf 6} at -0.3 6
\put {\bf 4} at 3.8 2.7
\put {\bf 7} at 5.8 7.7
}

\put {$4$} at 4.1 4.5
\put {$4$} at 4.45 3.3
\put {$3$} at 5.1 1.5
\put {$1$} at 6.45 3.3
\put {$7$} at 8.1 1.5
\put {$6$} at 8.1 4.5
\put {$5$} at 8.1 7
\put {$2$} at 8.95 8.3
\put {$7$} at 6.95 8.3
\put {$1$} at 5.95 6.3
\put {$5$} at 2.95 6.3
\put {$6$} at 0.95 6.3
\put {$1$} at 2.1 3.75
\put {$4$} at 0.95 1.8
\put {$5$} at 2.1 .75
\put {$7$} at 6.1 9

\setdots<2pt>

\plot 0 5 10 5
/

\endpicture
$$
\caption{Example of the Hammersley interacting fluid process}\label{fig:hamprocess}
\end{center}
\end{figure}

The Hammersley interacting fluid system $(M^{\nu}_t\,:\,t\geq 0)$ is a stochastic process with values in the space of positive, locally finite measures on $\R$. Its evolution is defined as follows: if there is a Poisson point with weight $\omega$ at a point $(x_0,t)$, then $M^\nu_{t}(\{x_0\}) = M^\nu_{t-}(\{x_0\}) + \omega$, and for $x> x_0$,
\begin{equation}\label{eq:evolM}
M^\nu_t((x_0,x]) = (M^\nu_{t-}((x_0,x]) - \omega)_+\,.
\end{equation}
Here, $M^\nu_{t-}$ is the ``mass distribution'' of the fluid at time $t$ if the Poisson point at $(x_0,t)$ would be removed. To the left of $x_0$ the measure does not change. In words, the Poisson point at $(x_0,t)$ moves a total mass $\omega$ to the left, to the point $x_0$, taking the mass from the first available fluid to the right of $x_0$. See Figure \ref{fig:hamprocess} below for a visualization, in case of atomic measures, of the process inside a space-time box. In this picture,  restricted to $[0,x]$, the measure $\nu$ consists of three atoms of weight $5$, $3$ and $7$. The measure $M^\nu_{t/2}$ consists of three atoms of weight $1$, $4$ and $6$, while at time $t$, it consists of one atom with weight $7$.

It is not true that the evolution $M_t$ is well defined for all measures $\nu$ (e.g. if we start with a finite number of particles to the left of $0$, every particle would be pulled instantaneously to $-\infty$). In this paper we follow the Aldous \& Diaconis \cite{AlDi} graphical representation in the last-passage model (compare to the result in the classical case, found in their paper):
\begin{prop}\label{prop:construction}
Let $\calN$ be the set of all positive, locally finite measures $\nu$ such that
\begin{equation}\label{eq:leftden}
\liminf_{y\to-\infty}\frac{\nu(y)}{y}>0\,.
\end{equation}
For each $\nu\in\calN$, the process defined by
\begin{equation}\label{eq:lastpassage}
L_{\nu}(x,t):=\sup_{z\leq x} \left\{ \nu(z) + L((z,0),(x,t))\right\}\ \ \ \ \ (x\in\R,\,\,t\geq 0)\,
\end{equation}
is well defined and the measure
\begin{equation}\label{eq:connection}
M_t^\nu((x,y]):=L_\nu(y,t)-L_\nu(x,t)\,,
\end{equation}
evolves according to the Hammersley interacting fluid system.
\end{prop}
\noindent{\bf Proof:} It follows immediately from the definition that $L_\nu(x,t)$ is increasing in $x$ and $t$, even if $L_\nu$ would not be finite everywhere. This implies that if we can prove that $L_\nu$ is finite on for example $\Z\times \Z_+$ with probability one, then almost surely, $L_\nu$ is finite everywhere. Therefore, we need only prove for any fixed point $(x,t)$ that $L_\nu(x,t)$ is finite with probability one. Use Theorem \ref{thm:shape} and \eqref{eq:sym} to see that
$$\E L\left((z,0),(x,t)\right)\sim \gamma \sqrt{t(x-z)}\ \ \ \ (z\to -\infty).$$
The Markov inequality and \eqref{eq:leftden} then give us the desired result.

We can show that $x\mapsto L_\nu((z,0),(x,t))$ is c\`adl\`ag using the fact that for $y\geq x$,
\[ L((z,0),(y,t)) \leq L((z,0),(x,t)) + L((x+,0),(y,t))\ \ \mbox{and}\ \ \lim_{y\downarrow x} L((x+,0),(y,t))=0,\]
where $x+$ means that you are not allowed to use a possible Poisson point directly above $x$. Therefore, $M^\nu_t$ is indeed a locally finite measure on $\R$. To see that $M^\nu_t$ follows the Hammersley interacting fluid dynamics that we have just defined, suppose that there is a Poisson point in $(x_0,t)$ with weight $\omega$. For $x<x_0$, this Poisson point has no effect, so $M^\nu_t = M^\nu_{t-}$ on $(-\infty,x_0)$. Clearly,
\[ M^\nu_{t}(\{x_0\}) = M^\nu_{t-}(\{x_0\}) + \omega.\]
If $x> x_0$, then the longest path to $(x,t)$ that attains the supremum in \eqref{eq:lastpassage} can either use the weight in $(x_0,t)$, which would give $L_\nu(x,t)=L_\nu(x_0,t)=L_\nu(x_0,t-)+\omega$, or it could ignore the weight in $(x_0,t)$, which would give $L_\nu(x,t)=L_\nu(x,t-)$. This proves that
\begin{eqnarray*}
M^\nu_{t}((x_0,x]) & = & L_\nu(x,t) - L_\nu(x_0,t)\\
& = & \max\left(L_\nu(x_0,t-)+\omega,L_\nu(x,t-)\right) - L_\nu(x_0,t-) - \omega\\
& = & \left(M_{t-}^\nu((x_0,x])-\omega\right)_+.
\end{eqnarray*}
\hfill$\Box$\\

Assume that we have a probability measure defined on $\calN$ and consider $\nu\in \calN$ as a realization of this probability measure. We say that $\nu$ is time invariant for the Hammersley interacting fluid process (in law) if
$$M^\nu_t\stackrel{\cal D}{=}M^\nu_0=\nu\,\,\,\mbox{ for all }\,\,\,t\geq 0\,.$$
In this case, we also say that the underlying probability measure on $\calN$ is an equilibrium measure. Let $\alpha \in (\pi,3\pi/2)$ and define the measure in $\calN$
$$\nu_\alpha((x,y])=B_\alpha\left((x,0),(y,0)\right)\mbox{ for }\,\,\, x\leq y\in\R\,.$$
Compare this to \eqref{eq:defnua} ($\nu_\alpha = \nu^\alpha_0$).
It was this interplay between the longest path description and the equilibrium interacting particle system that proved very fruitful in the results for the classical Hammersley process in Cator \& Groeneboom \cite{CaGr2}. We will attempt the same in the interacting fluid system, but since the equilibrium solution is not explicitly known, we needed to develop new tools and ideas, which in fact also had interesting applications for the classical case. Of course, as an immediate consequence of Proposition \ref{prop:BusMarkov}, we do have the following
\begin{coro}\label{cor:equilibrium}
The random measures $\nu_\alpha = \nu_0^{\alpha}$ \eqref{eq:defnua} are all equilibrium measures for the Hammersley interacting fluid process.
\end{coro}

\section{Ergodicity and uniqueness of the equilibrium measure}\label{sec:ergodic}
Let $\alpha \in (\pi,3\pi/2)$ and define the measure in $\calN$
$$\nu_\alpha((x,y])=B_\alpha\left((x,0),(y,0)\right)\mbox{ for }\,\,\, x\leq y\in\R\,.$$
Compare this to \eqref{eq:defnua} ($\nu_\alpha = \nu^\alpha_0$). To prove the next theorem about time invariance, we need to define the following exit points:
\begin{equation}\label{def:Zalpha}
Z_\alpha(x,t) \mbox{ is the crossing point of }\varpi_\alpha(x,t)\ \mbox{and } \R\times\{0\}.
\end{equation}
This was already used in the proof of Proposition \ref{prop:BusMarkov}. Analogously,
\begin{equation}\label{def:Zalpha*}
Z^*_\alpha(x,t) \mbox{ is the crossing point of }\varpi_\alpha(x,t)\ \mbox{and } \{0\}\times\R.
\end{equation}
Finally, if $\nu\in \calN$, we define
\begin{equation}\label{def:Znu}
Z_\nu(x,t) = \sup\{z\leq x\ :\ L_\nu(x,t) = \nu(z) + L((z,0),(x,t))\}.
\end{equation}
This says that $Z_\nu(x,t)$ is the right-most point where the supremum in the definition of $L_\nu(x,t)$ is attained. This definition is slightly subtle. The fact that the supremum in the definition of $L_\nu(x,t)$ is attained, relies on the fact that on a compact interval, the sum of a non-decreasing right-continuous function ($\nu(z)$) and a non-increasing {\em left}-continuous function ($z\mapsto L((z,0),(x,t))$) attains its maximum. Now we have to show that for all $(y,s)$, the supremum over $z\leq y$ can actually be restricted to a compact set. For a countable set of $(x,t)$'s this can be done as in the proof of Proposition \ref{prop:construction}, using \eqref{eq:leftden} and Theorem \ref{thm:shape} for $L$. Then we can conclude the desired compactness property for all $(y,s)$ by using the inequality
\[ \nu(z) + L((z,0),(y,s))\leq \nu(Z_\nu(x,t)) + L((Z_\nu(x,t),0),(y,s))\ \ \forall\ z\leq Z_\nu(x,t), y\geq x, s\leq t.\]

We have the following proposition.
\begin{prop}\label{prop:exitpoint}
For all $x\in \R$ and $t\geq 0$, we have
\[ Z_\alpha(x,t) = Z_{\nu_\alpha}(x,t).\]
\end{prop}
\noindent{\bf Proof:} From \eqref{eq:zZalpha} we immediately get that $Z_{\nu_\alpha}(x,t)\geq Z_\alpha(x,t)$. Now suppose $Z_{\nu_\alpha}(x,t)> Z_\alpha(x,t)$. Since $L(Z_\alpha(x,t),(x,t)) = L(Z_{\nu_\alpha}(x,t),(x,t)) + B_\alpha(Z_{\nu_\alpha}(x,t), Z_\alpha(x,t))$, we could have chosen the $\alpha$-ray starting at $(x,t)$ through the point $Z_{\nu_\alpha}(x,t)$, which would make it strictly lower than $\varpi_\alpha(x,t)$, going through $Z_\alpha(x,t)$.\hfill $\Box$\\

Since the $\alpha$-ray $\varpi_\alpha(\0)$ has asymptotic direction $(\cos(\alpha),\sin(\alpha))$ with probability 1, we can easily see that
\begin{equation}\label{eq:vanishexitpoint}
\lim_{t\to\infty}\P\big(|Z_{\alpha}(t,(\tan\alpha) t)|\geq\epsilon t\big)=\lim_{t\to\infty}\P\big(|Z^*_{\alpha}(t,(\tan\alpha) t)|\geq\epsilon t\big)=0\,.
\end{equation}
However, we can also control $Z_\nu$ for more general $\nu$. Compare the following lemma to Lemma 3.3 in Ferrari, Martin \& Pimentel \cite{FeMaPi}.
\begin{lem}\label{lem:genexitcontrol}
Suppose $\nu\in \calN$. Assume that
\begin{equation}\label{eq:genexitcontrol}
\liminf_{z\to-\infty}\frac{\nu((z,0))}{-z}\geq\frac{\gamma}{2}\sqrt{\tan\alpha}\,\,\mbox{ and that }\,\,\limsup_{z\to\infty}\frac{\nu([0,z])}{z}\leq\frac{\gamma}{2}\sqrt{\tan\alpha}\,.
\end{equation}
Then, with probability one,
$$\lim_{t\to\infty}\frac{Z_{\nu}\big(t,(\tan\alpha) t\big)}{t}=0\,.$$
\end{lem}
\noindent{\bf Proof:}
Using the transformation \eqref{eq:sym}, we can assume, without loss of generality, that $\alpha=5\pi/4$. The proof of this lemma is based on the following estimate for the shape function $f(x,t)=\gamma\sqrt{xt}$: for $s\in[0,t]$
\begin{equation}\label{eq:curvature}
f(t,t)-f(t-s,t)\geq \frac{\gamma}{2}s + \frac{\gamma}{8}\frac{s^2}{t}\,.
\end{equation}
Now, fix $\eps>0$ and suppose that $Z_{\nu}(t):=Z_{\nu}(t,t)\geq \eps t + 1$. Then there exists $k\in[\eps t,t]\cap {\mathbb N}$ and $n=\lfloor t\rfloor$ such that
$$\nu(k+1)+L\big((k,0),(n+1,n+1)\big)-L\big((0,0),(n,n)\big)\geq 0\,.$$
By adding $f(n,n)-f(n-k,n)-\frac{\gamma}{2}k$ to both sides of the above inequality and applying (\ref{eq:curvature}), we get
\begin{eqnarray}
\nonumber\nu(k+1)-\frac{\gamma}{2}k\ + & &\\
\nonumber L\big((k,0),(n+1,n+1)\big)-f(n-k,n)\ + & &\\
\nonumber f(n,n)-L\big((0,0),(n,n)\big)\ &\geq&\,\,\frac{\gamma}{8}\frac{k^2}{n}\\
\label{eq:curv}& \geq & \frac{\gamma}{8}\eps^2n\,.
\end{eqnarray}
By \eqref{eq:genexitcontrol}, we know that
\[ \limsup_{t\to \infty} \sup_{s\in[\eps t,t]} \frac{\nu(s+1)-\frac{\gamma}{2}s}{t} \leq 0\,.\]
Therefore, by \eqref{eq:curv}, if the set $\{t\geq 0\,:\, Z_{\nu}(t,t)\geq \eps t+1\}$ is unbounded, it follows that for some small $\eta>0$, the events
\[\{|L\big((k,0),(n+1,n+1)\big)-f(n-k,n)|> \eta n\}\]
happen infinitely often for $n\geq 1$ and $0\leq k\leq n$. Using Borel-Cantelli, this will have zero probability if for all $\eta>0$,
\begin{equation}\label{eq:borcant}
\sum_{n=1}^\infty\, \sum_{k=0}^n\, \P\left(|L\big((k,0),(n+1,n+1)\big)-f(n-k,n)|> \eta n\right) < +\infty.
\end{equation}
Theorem \ref{thm:shape} gives us some control on the fluctuations of $L$ about its asymptotic shape. Note that for $n$ big enough, and $0\leq k\leq n$, we have
\[\begin{array}{ll}
\P\big(|L\big((k,0),(n+1,n+1)\big)\hspace{-0.2cm}&-\ f(n-k,n)|>\eta n\big) = \vspace{0.3cm}\\
& = \P\left(|L\big(\0,\sqrt{(n+1-k)(n+1)}(1,1)\big) - f(n-k,n)|>\eta n\right).
\end{array}\]
Define $r=\sqrt{(n+1-k)(n+1)}$ and $u=n^{2/3}$. If we choose $n$ large enough, we can make sure that $u<\eta n/2$ and
\[|f(n-k,n) - \gamma r|<\eta n/2.\]
Also, for $n$ large enough, we have that $u\in \left[\,c_3 \sqrt{r}\log^2 r\,,\,c_4 r^{3/2}\log r\,\right]$ (see Theorem \ref{thm:shape}). This implies, using Theorem \ref{thm:shape}, that there exist $c_1,c_2>0$ such that for $n$ large enough and $0\leq k\leq n$,
\begin{eqnarray*}
 \P\big(|L\big((k,0),(n+1,n+1)\big) - f(n-k,n)|>\eta n\big)& \leq & \P\big(|L\big(\0,(r,r)\big) - \gamma r|>u\big)\\
 &\leq & c_1\exp\left(-c_2n^{1/6}/\log(n)\right).
\end{eqnarray*}
This clearly proves \eqref{eq:borcant}. The proof that the set $\{t\geq 0\,:\, Z_{\nu}(t,t)\leq -\eps t - 1\}$ is bounded with probability 1 follows the same line.\hfill $\Box$\\

Now we can prove the most important result of this Section.
\begin{thm}\label{thm:invariant}
If we start the Hammersley interacting fluid system with $\nu_\alpha$ then
\[M_t^{\nu_\alpha}\stackrel{\cal D}{=} \nu_\alpha\,\,\,\mbox{ for all }\,\,\,t\geq 0\,.\]
The process $x\mapsto \nu_\alpha(x)$ is stationary and ergodic and its intensity is given by
\begin{equation}\label{eq:gamma}
\E \nu_\alpha(1)=\frac{\gamma(F)}{2}\sqrt{\tan\alpha}\,.
\end{equation}
Finally, consider a random $\nu\in \calN$, which is time invariant, and which defines a stationary and ergodic process on $\R$. Define $\alpha\in (\pi,3\pi/2)$ by
\[ \alpha =\arctan \left(\frac{2}{\gamma(F)}\E \nu(1)\right)^2.\]
Then $\nu \stackrel{\cal D}{=}  \nu_\alpha$.
\end{thm}
\noindent {\bf Proof:} The first statement is an immediate consequence of Proposition \ref{prop:construction} and Proposition \ref{prop:BusMarkov}. The fact that $x\mapsto \nu_\alpha(x)$ is stationary and that $0<\E(\nu_\alpha(1))<+\infty$ follows directly from Proposition \ref{prop:Busemann}. Ergodicity will follow from Proposition \ref{prop:mixing}, but won't be needed for the rest of this proof. Fix $\alpha\in(\pi,3\pi/2)$ and set $\rho=\rho(\alpha):=\sqrt{\tan\alpha}$. Since the model is invariant under the map $(x,t)\to(\rho x,t/\rho)$,
$$\E \nu_\alpha(1)= \E \nu _{5\pi/4}(\rho) =\E \nu _{5\pi/4}(1)\sqrt{\tan\alpha}\,.$$
Now, for all $t\geq 0$,
$$L_{5\pi/4}(t,t)=\nu_{5\pi/4}\left(Z_{\nu_{5\pi/4}}(t,t)\right)+L\big((Z_{\nu_{5\pi/4}}(t,t),0),(t,t)\big)\,.$$
By Proposition \ref{prop:exitpoint}, $Z_{\nu_{5\pi/4}}(t,t)=Z_{5\pi/4}(t,t)=:Z_{5\pi/4}(t)$. If $Z_{5\pi/4}(t)\geq 0$ then
$$0\leq L_{\nu_{5\pi/4}}(t,t)-L\big((0,0),(t,t)\big)\leq \nu_{5\pi/4}\left(Z_{5\pi/4}(t)\right)\,.$$
On the other hand, if $Z_{5\pi/4}(t)< 0$ then $Z^*_{5\pi/4}(t)\geq 0$. Now define
\begin{equation}\label{eq:defnu*}
\nu^*_\alpha(x) = B_\alpha(\0,(0,x)).
\end{equation}
From the additivity of the Busemann function, we know that
$$L_{5\pi/4}(t,t)=\nu^*_{5\pi/4}\left(Z^*_{5\pi/4}(t,t)\right)+L\big((0,Z^*_{{5\pi/4}}(t,t)),(t,t)\big)\,.$$
Finally we obtain that, denoting  $Z^*_{5\pi/4}(t,t)=:Z^*_{5\pi/4}(t)$,
$$0\leq L_{\nu_{5\pi/4}}(t,t)-L\big((0,0),(t,t)\big)\leq \max\left\{\nu_{5\pi/4}(Z_{5\pi/4}(t)),\nu^*_{5\pi/4}(Z^*_{5\pi/4}(t))\right\}\,.$$
Clearly, from Proposition \ref{prop:Busemann}\eqref{prop:Busesym}  (symmetry) it follows that $\nu^*_{5\pi/4} \stackrel{\cal D}{=} \nu_{5\pi/4}$. Together with \eqref{eq:vanishexitpoint}, this yields that
\begin{equation}\label{eq:int}
\frac{|L_{\nu_{5\pi/4}}(t,t)-L\big((0,0),(t,t)\big)|}{t} \stackrel{\cal D}{\longrightarrow} 0 \,.
\end{equation}
Here we use the fact that for any $\eta>0$, with probability one, $\max\left\{Z_{5\pi/4}(t),Z^*_{5\pi/4}(t)\right\}\leq \eta t$ for $t$ large enough. Then we can use stationarity and the fact that $\E \nu_\alpha(1)<+\infty$ to conclude that for any $\eps>0$, there exists $\eta>0$, such that $\P(\nu_{5\pi/4}(\eta t)\geq \eps t)\leq \eps$.

By Corollary \ref{cor:equilibrium}, we get
\begin{eqnarray*}
\E L_{\nu_{5\pi/4}}(t,t) & = & \E(L_{\nu_{5\pi/4}}(t,t) - L_{\nu_{5\pi/4}}(0,t)) + \E L_{\nu_{5\pi/4}}(0,t)\\
& = & \E(\nu_{5\pi/4}(t)) + \E(\nu^*_{5\pi/4}(t))\\
& = & 2t\E(\nu_{5\pi/4}(1)).
\end{eqnarray*}
The sub-additive ergodic theorem applied to $L_{5\pi/4}(t,t)$ implies that, with probability one,
$$\lim_{t\to\infty}\frac{L_{5\pi/4}(t,t)}{t}= 2\E \nu_{5\pi/4}(1)\,.$$
Combining this with Theorem \ref{thm:shape} and (\ref{eq:int}), one gets (\ref{eq:gamma}).\\

Now we need to address the uniqueness of $\nu_\alpha$. Suppose $\nu\in\calN$ is ergodic and time invariant. Define $Z(t)=Z_\nu(t,t\tan(\alpha))$ and $Z_h(t)=Z_\nu(t+h,t\tan(\alpha))$. Now define
\[ \tilde{Z}(t) = \argmax_{z\leq t}\left( \nu(z) + L((-t+z,-t\tan(\alpha)),\0)\right)\]
and
\[ \tilde{Z}_h(t) = \argmax_{z\leq t+h}\left( \nu(z) + L((-t+z,-t\tan(\alpha)),(h,0))\right).\]
Here, we take the right-most location of the maximum. The intuition for $\tilde{Z}(t)$ and $\tilde{Z}_h(t)$ is that we place the origin at $(-t,-t\tan(\alpha))$, and look at the exit-point for the path that starts at $(-t,-t\tan(\alpha))$, picks up mass from $\nu$ and then goes to $\0$, resp. $(h,0)$. Clearly, we have
\[ (Z(t), Z_h(t))\stackrel{\cal D}{=} (\tilde{Z}(t), \tilde{Z}_h(t)).\]
Since $\nu$ is ergodic, and by our choice of $\alpha$, $\nu$ satisfies \eqref{eq:genexitcontrol}. Since the proof of Lemma \ref{lem:genexitcontrol} uses a Borel-Cantelli type argument, it is not hard to see that we can use the same ideas to prove
\[ (\tilde{Z}(t),\tilde{Z}_h(t))/t \stackrel{\rm a.s.}{\longrightarrow} (0,0).\]
This means that the two paths $\varpi((-t+\tilde{Z}(t),-t\tan(\alpha)),\0)$ and $\varpi((-t+\tilde{Z}_h(t),-t\tan(\alpha)),(h,0))$ will converge in any bounded box to the $\alpha$-rays $\varpi_\alpha(\0)$ and $\varpi_\alpha((h,0))$ respectively (this follows from Theorem \ref{thm:existray}(2)). However, these two $\alpha$-rays will coalesce, which means that with probability 1, there exists $t_0>0$ such that for all $t\geq t_0$, the two converging paths coalesce, which in turn implies that $\tilde{Z}(t)=\tilde{Z}_h(t)$ (because they are both the right-most point where the maximum takes place and, as soon as they coalesce, they get the same exit point). Now define
\[ \tilde{L}(t) = \sup_{z\leq t}\left( \nu(z) + L((-t+z,-t\tan(\alpha)),\0)\right)\]
and
\[ \tilde{L}_h(t) = \sup_{z\leq t+h}\left( \nu(z) + L((-t+z,-t\tan(\alpha)),(h,0))\right).\]
We also have that
\[ (\tilde{L}(t), \tilde{L}_h(t)) \stackrel{\cal D}{=} (L_\nu(t,t\tan(\alpha)),L_\nu(t+h,t\tan(\alpha))).\]
Furthermore, if $t\geq t_0$, then
\[ \tilde{L}_h(t) - \tilde{L}(t) = B_\alpha(\0,(0,h))=\nu_\alpha((0,h]).\]
This proves that
\begin{eqnarray*}
M^\nu_{t\tan(\alpha)}((t,t+h]) & = & L_\nu(t+h,t\tan(\alpha)) - L_\nu(t,t\tan(\alpha))\\
 & \stackrel{\cal D}{\longrightarrow} & \nu_\alpha((0,h]).
\end{eqnarray*}
Since $\nu$ is time invariant and ergodic, we see that
\[ \nu((0,h]) \stackrel{\cal D}{=} M^\nu_{t\tan(\alpha)}((t,t+h]) \stackrel{\cal D}{=} \nu_\alpha((0,h]).\]
In principle, we need to show convergence for a finite number of $h$'s simultaneously, but it is not hard to see that the ideas we used can be extended to that case, at the cost of some notational burden.  Note that we have proved that for any deterministic $\nu$ satisfying (\ref{eq:genexitcontrol}), $M^\nu_{(\tan\alpha)t}([t,t+h])$ converges in distribution to $\nu_{\alpha}(h)$, as a process in $h$. This shows that in a rarefaction fan, the fluid process converges locally to the correct equilibrium process (local equilibrium).
\mbox{}\hfill $\Box$\\

\begin{coro}\label{cor:nu*}
\[ \E(\nu_\alpha^*(1)) = \frac{\gamma(F)}{2\sqrt{\tan(\alpha)}}.\]
In particular, for all $\alpha\in (\pi,3\pi/2)$, we have
\[ \E(\nu_\alpha(1))\cdot \E(\nu_\alpha^*(1)) = \frac{\gamma(F)^2}{4}.\]
\end{coro}
\noindent {\bf Proof:} Remember that $\nu^*_\alpha(x) = B_\alpha(\0,(0,x))$. For $\alpha = 5\pi/4$, the result follows Proposition \ref{prop:Busemann}\eqref{prop:Busesym}  (symmetry) and Theorem \ref{thm:invariant}. Now use the map $(x,t)\mapsto (\rho x,t/\rho)$ to see that
\[ \nu^*_\alpha(x) \stackrel{\cal D}{=} \nu^*_{5\pi/4}\left(x/\sqrt{\tan(\alpha)}\right).\]
\mbox{}\hfill $\Box$\\

For the classical Hammersley model, we know that if $\bar\nu_\lambda$ is a Poisson counting process of intensity $\lambda$, then $\bar\nu_\lambda$ is time invariant and ergodic. Therefore it must be equal in distribution to $\nu_\alpha$, for some $\alpha \in (\pi, 3\pi/2)$. We also know that in the classical Hammersley process,
\[ L_{\bar{\nu}_\lambda}(0,t) \stackrel{\cal D}{=} \bar{\nu}_{1/\lambda}((0,t]).\]
So for $\lambda=1$, we get
\[  L_{\bar{\nu}_1}(0,t) \stackrel{\cal D}{=} \bar{\nu}_1((0,t]).\]
Since we know that for any $\alpha \in (\pi, 3\pi/2)$,
\[ L_{\nu_\alpha}(0,t) = B_\alpha(\0,(0,t))\stackrel{\cal D}{=}B_{5\pi/2-\alpha}(\0,(t,0)),\]
using \ref{prop:Busemann}\eqref{prop:Busesym}, we can now use Corollary \ref{cor:nu*} to conclude that $\nu_{5\pi/4} = \bar{\nu}_1$.
Consequently,
$$1=\E \bar\nu_1(1)=\E \nu_{5\pi/4}(1)=\frac{\gamma(\delta_1)}{2}\sqrt{\tan(5\pi/4)}=\frac{\gamma(\delta_1)}{2}\,,$$
which proves that $\gamma(\delta_1)=2$. We remark that the proof that the Poisson process is time invariant does not depend on the value of $\gamma(\delta_1)$. It only relies on an explicit calculation of the generator associated to  $M_t$. (See Lemma 8 of Aldous \& Diaconis \cite{AlDi}, or  Theorem 3.1 of Cator \& Groeneboom \cite{CaGr1}.)
\begin{coro}\label{coro:original}
In the classical Hammersley model, we have that $\gamma(1)=2$ and that $\nu_\alpha\stackrel{\cal D}{=} \bar\nu_{\lambda(\alpha)}$ where $\lambda(\alpha)=\sqrt{\tan\alpha}$.
\end{coro}

For general weight distributions $F$, we were not able to get more information on $\nu_\alpha$ (not even a guess for a good candidate). In particular, we do not know how to calculate $\gamma(F)$. This does seem to be the most important contribution of the interacting fluid representation: once we have a good candidate for $\nu_\alpha$, we can check it by showing that it is invariant under the evolution of the interacting fluid. In fact, even in the results for the classical Hammersley case found in Aldous \& Diaconis \cite{AlDi} and Cator \& Groeneboom \cite{CaGr1,CaGr2} , this is where the interacting particle process proves its worth.

\subsection{Mixing property of $\nu_\alpha$}

We will show that the measure $\nu_\alpha$ has the following mixing property, usually called strong mixing in dynamical systems. We consider the $\sigma$-algebra $\calF = \sigma\{ \nu_\alpha((a,b])\ :\ a\leq b \in \R\}$ on the sample space $\Omega$, defined by the compound Poisson process. We can define the translation $\tau_t$ as an $\calF$-measurable map from $\Omega$ to $\Omega$, simply by translating all Poisson points by the vector $(t,0)$.
\begin{prop}\label{prop:mixing}
For each $\alpha\in(\pi,3\pi/2)$, $\nu_\alpha$ satisfies
\begin{equation}\label{eq:defmixing}
 \forall\ A,B\in \calF:\ \lim_{t\to \infty}\ \P(A\cap \tau^{-1}_t(B)) = \P(A)\P(B).
 \end{equation}
In particular, this implies that $\nu_\alpha$ is ergodic.
\end{prop}
\noindent{\bf Proof:} From translation invariance and a standard approximation of sets in $\calF$, it is enough to prove \eqref{eq:defmixing} for all $A,B\in \calF_h:=\sigma\{ \nu_\alpha((a,b])\ :\ a\leq b \in [0,h]\}$. Consider the paths $\varpi((-t,-t\tan(\alpha)),\0)$ and $\varpi((-t,-t\tan(\alpha)),(h,0))$. Almost surely, these paths will converge to $\varpi_\alpha(\0)$ and $\varpi_\alpha((h,0))$, respectively, on any finite box. This means, that if we define for $a,b\in [0,h]$
\[ \nu^{(t)}_\alpha((a,b]) = L((b,0),(-t,-t\tan(\alpha))) - L((a,0),(-t,-t\tan(\alpha))),\]
then for $t$ big enough, we have $\nu^{(t)}_\alpha = \nu_\alpha|_{[0,h]}$. Clearly, $\tau^{-1}_{t+h}(B)$ is independent of $\nu_\alpha^{(t)}$, since they depend on the Poisson process to the left respectively to the right of the line $\{-t\}\times \R$. Define the event
\[ C_t = \{ \forall\ s\geq t:\ \nu^{(s)}_\alpha = \nu_\alpha|_{[0,h]}\}\]
and denote $A^{(t)}$ the counterpart of the event $A$ in $\calF_h^{(t)}:=\sigma\{ \nu^{(t)}_\alpha((a,b])\ :\ a\leq b \in [0,h]\}$; although it is intuitively clear what is meant, we will make this more precise at the end of the proof. Then
\begin{eqnarray*}
|\P(A\cap \tau_{t+h}^{-1}(B)) - \P(A)\P(B)| & \leq & |\P(A\cap \tau_{t+h}^{-1}(B)\cap C_t) - \P(A)\P(B)| + \P(C^c_t)\\
& = & |\P(A^{(t)}\cap \tau_{t+h}^{-1}(B)\cap C_t) - \P(A)\P(B)| + \P(C^c_t)\\
& \leq & |\P(A^{(t)}\cap \tau_{t+h}^{-1}(B)) - \P(A)\P(B)| + 2\P(C^c_t)\\
& = & \P(B)|\P(A^{(t)}) - \P(A)| + 2\P(C^c_t)\\
& \leq & 4\P(C^c_t).
\end{eqnarray*}
The proposition now follows from the fact that $\P(C^c_t)\to 0$\footnote{We note that, in the classical model, we have independent increments even if the probability of the event $C_t^c$ does not decay to $0$ very fast. This indicates that, to show mixing by using these events may not be the best strategy.}.

To see what is meant by $A^{(t)}$, we define the index-set $I=\{ \nu_\alpha((a,b])\ :\ a\leq b \in [0,h]\}$ and $I^{(t)}=\{ \nu^{(t)}_\alpha((a,b])\ :\ a\leq b \in [0,h]\}$. There is a canonical bijection $i:I\to I^{(t)}$. Define $\cal B$ as the product $\sigma$-algebra on $\R^I$, and likewise ${\cal B}^{(t)}$. Extend the canonical map $i$ such that $i:\R^I\to\R^{I^{(t)}}$. Define the map
\[ \phi: \Omega \to \R^I : \omega \mapsto \{\nu_\alpha((a,b])(\omega)\ :\ a\leq b \in [0,h]\},\]
and likewise $\phi_t:\Omega \to \R^{I^{(t)}}$. We know that $\calF_h = \phi^{-1}(\cal B)$ and $\calF_h^{(t)}=\phi_t^{-1}(\cal B^{(t)})$. This means that there exists $U\in \cal B$, such that $A=\phi^{-1}(U)$. We define $A^{(t)} = \phi_t^{-1}(i(U))$.
\hfill $\Box$\\

\section{Central limit theorems for the Busemann function in the classical model}\label{sec:scalBus}
Our geometrical approach yields a very explicit description of the fluctuations of the Busemann function in the classical Hammersley model. We first notice the following relations that will be derived from Proposition \ref{prop:Busemann} and Proposition \ref{prop:BusMarkov}:
\begin{prop}\label{PropCLT}
Consider the Hammersley last-passage model with random weights and recall $\vec\beta:=(\cos\beta,\sin\alpha)$.
\begin{itemize}
\item If $\beta\in[0,\pi/2]$ then, as processes,
$$B_\alpha(\0,\cdot\vec\beta)= L_{\nu_\alpha}(\cdot\vec\beta)\,;$$
\item If $\beta\in[\pi/2,\pi]$ then for each $t\geq 0$
$$B_\alpha(\0,t\vec\beta)\stackrel{\cal D}{=} \nu^*_\alpha\left(t\sin\beta\right)-\nu_\alpha\left(-t\cos\beta\right)\,.$$
\end{itemize}
\end{prop}
\noindent{\bf Proof:} The first statement follows directly from Proposition \ref{prop:BusMarkov}. To obtain the second relation note that, by additivity and anti-symmetry (Proposition \ref{prop:Busemann}),
\begin{eqnarray}
\nonumber B_\alpha\big(\0,t\vec\beta\big)&=&B_\alpha\big(\0,(t\cos\beta,0)\big)+B_\alpha\big((t\cos\beta,0),t\vec\beta\big)\\
\nonumber &=&B_\alpha\big((t\cos\beta,0),t\vec\beta\big)-B_\alpha\big((t\cos\beta,0),\0\big)\,.
\end{eqnarray}
Now use translation invariance, translating all four points by $-(t\cos\beta,0)$.\hfill $\Box$\\

In the classical set up, Corollary \ref{coro:original}, together with Proposition \ref{PropCLT}, implies that:
\begin{itemize}

\item[(1)] if $\beta=0$ then
$$B_\alpha\big(\0,t\vec\beta\big)\stackrel{\cal D}{=}X_\alpha(t\cos\beta)\,,$$
where $X_\alpha$ is a Poisson process of intensity $\lambda(\alpha)$.
\item[(2)] if $\beta\in(0,\pi/2)$ then
$$B_\alpha\big(\0,\cdot\vec\beta\big)\stackrel{\cal D}{=}L_{\bar\nu_{\lambda(\alpha)}}(\cdot\vec\beta)\,,$$
where $\lambda=\lambda(\alpha)=\sqrt{\tan\alpha}$;
\item[(3)] if $\beta\in[\pi/2,\pi]$ then for each $t\geq 0$
$$B_\alpha\big(\0,t\vec\beta\big)\stackrel{\cal D}{=}Y_\alpha(t\sin\beta)-X_\alpha(-t\cos\beta)\,,$$
where $Y_\alpha$ and $X_\alpha$ are two independent one dimensional Poisson processes variables of intensity $1/\lambda(\alpha)$ and $\lambda(\alpha)$, respectively.
\end{itemize}

Baik and Rains (2001) proved the following central limit theorem for the Hammersley classical model with external ``sources''. Let $\Phi(x)$ be the standard normal distribution function, and let $F_0(x)$ be the zero mean Tracy-Widon type distribution function introduced in Definition 2 of \cite{BaRa}.
\begin{itemize}
\item[(4)] If $\lambda\in(0,1)$ then
$$\lim_{s\to\infty}\P\left(\frac{L_{\bar\nu_{\lambda}}(s,s)-(1/\lambda+\lambda)s }{(\sqrt{1/\lambda-\lambda})s^{1/2}}\leq x\right)= \Phi(x)\,;$$
\item[(5)] If $\lambda=1$ then
$$\lim_{s\to\infty}\P\left(\frac{L_{\bar\nu_{\lambda}}(s,s)-2s }{s^{1/3}}\leq x\right)= F_0(x)\,;$$
\item[(6)] If $\lambda>1$ then
$$\lim_{s\to\infty}\P\left(\frac{L_{\bar\nu_{\lambda}}(s,s)-(\lambda+1/\lambda)s }{(\sqrt{\lambda-1/\lambda})s^{1/2}}\leq x\right)= \Phi(x)\,.$$
\end{itemize}

These results naturally lead us to a central limit theorem for the Busemann function
$$B(\beta,\cdot):=B_{5\pi/4}\big(\0,(\cdot\cos\beta,\cdot\sin\beta)\big)\,.$$
By \eqref{eq:sym}, w.l.o.g. we can restrict ourselves to $\alpha=5\pi/4$. In this case we have
$$\E B(\beta,t)=(\cos\beta+\sin\beta)t\,.$$
\begin{coro}\label{coro:TCLBus}
Consider the Hammersley classical last-passage model. Then
$$\lim_{t\to\infty}\P\left(\frac{B\big(\pi/4,t\big)-\sqrt{2}t }{2^{-1/6}t^{1/3}}\leq x\right)= F_0(x)\,,$$
while for $\beta\in[0,\pi/4)$
$$\lim_{t\to\infty}\P\left(\frac{B\big(\beta,t\big)-(\cos\beta+\sin\beta)t }{(\sqrt{\cos\beta-\sin\beta})t^{1/2}}\leq x\right)= \Phi(x)\,,$$
and for $\beta\in(\pi/4,\pi]$
$$\lim_{t\to\infty}\P\left(\frac{B\big(\beta,t\big)-(\cos\beta+\sin\beta)t }{(\sqrt{\sin\beta-\cos\beta})t^{1/2}}\leq x\right)= \Phi(x)\,.$$
\end{coro}
\noindent{\bf Proof:} The statement for $\beta=0$ and $\beta\in(\pi/2,\pi]$ follows (1) and (3), together with the central limit theorem. Now, let us take $\beta\in(0,\pi/2)$. Consider the map \eqref{eq:sym} with $\rho(\beta)=\sqrt{\tan{\beta}}$. Then $B(\beta,t)\stackrel{\cal D}{=}L_{\bar\nu_\rho}(s,s)$, where $s=(\sqrt{\sin\beta\cos\beta})t$, and for $\beta\in(0,\pi/4]$ we have $\rho\in(0,1]$ while for $\beta\in(\pi/4,\pi]$ we have $\rho\in(1,\infty)$.\hfill $\Box$\\

\subsection{The crossing formula}
In the classical model, the exit point formula for the equilibrium regime, proved by Cator and Grooeneboom \cite{CaGr2}, is
\begin{equation}\label{eq:exit}
{\rm Var} L_{\bar{\nu}_\lambda}(x,t)=-\lambda x + \frac{t}{\lambda}+2\lambda\E Z_{\bar{\nu}_\lambda}(x,t)_+\,,
\end{equation}
where ${\rm Var} X$ is the variance of $X$, and $X_+:=\max\{X,0\}$. Together with (2) and Proposition \ref{prop:exitpoint}, this relates the variance of the Busemann function at $(x,t)$ with the position of the crossing point of the $\alpha$-ray starting at $(x,t)$:
\begin{equation}\label{eq:cross}
{\rm Var} B_\alpha(\0,(x,t))=-(\sqrt{\tan\alpha}) x + \frac{t}{\sqrt{\tan\alpha}}+2(\sqrt{\tan\alpha})\E Z_\alpha(x,t)_+\,.
\end{equation}
In particular,
$${\rm Var} B_{5\pi/4}(\0,(t,t))=2\E Z_{5\pi/4}(t,t)_+\,.$$
We note that this crossing point formula can be seen as a version of the scaling identity $\xi=2\chi$, where $\chi$ and $\xi$ are the critical exponents that measure the order of magnitude of the fluctuations of Busemann functions and crossing points, respectively.

\section{The multi-class process and second class particles}\label{sec:multiclass}

For two positive measures $\nu$ and $\bar\nu$ on $\R$, we say that $\bar\nu$ dominates $\nu$, notation $\bar\nu\geq\nu$, whenever $\bar\nu(I)\geq \nu(I)$ for all measurable $I\subseteq \R$.
\begin{prop}\label{prop:coupling}
Suppose we have two measures $\nu, \bar\nu\in {\cal N}$ such that $\bar\nu\geq \nu$. Define the corresponding interacting fluid system as $M^\nu_t$ and $M^{\bar\nu}_t$, using the same weighted Poisson process (basic coupling). Then $M^{\bar\nu}_t\geq M^\nu_t$ (as measures). If $\bar\nu(z)=\nu(z)$ for all $z<0$, then  $M^{\bar\nu}_t([0,x])-M^\nu_t([0,x])$ is non-increasing in $t$ for all $x\geq 0$.
\end{prop}
\noindent{\bf Proof:} Fix an interval $[-K,K]$ and a time $t$. There exists (a random) $M>0$ such that $M^\nu_t$ and $M^{\bar\nu}_t$ restricted to $[-K,K]$ only depend on Poisson points in $[-M,K]\times [0,t]$ and on $\nu$ and $\bar\nu$ restricted to $[-M,K]$ (it is not hard to see that we can take $M=Z_\nu(-K,t)$). This means that we are only dealing with a finite number of Poisson points, so if we can prove that the premise ``$M^{\bar\nu}_s\geq M^\nu_s$ for all $s<t$'' implies that $M^{\bar\nu}_t\geq M^\nu_t$, we will have proved the first statement, since it is obviously true for $t=0$. Suppose there exists a Poisson point at $(x_0,t)$ with weight $\omega$ for some $x_0\in [-M,K]$, since otherwise the implication is immediate. We then know, using Proposition \ref{eq:leftden} and \eqref{eq:evolM}, that if $x_0< x\leq y$,
\begin{eqnarray*}
M^{\nu}_t((x,y]) & = & (M^{\nu}_{t-}((x_0,y]) - \omega)_+  - (M^{\nu}_{t-}((x_0,x]) - \omega)_+ \\
& \leq & (M^{\bar\nu}_{t-}((x_0,y]) - \omega)_+  - (M^{\bar\nu}_{t-}((x_0,x]) - \omega)_+ \\
& = & M^{\bar\nu}_t((x,y]).
\end{eqnarray*}
The inequality follows from the fact that if $A\geq B$ and $\tilde{A}\geq \tilde{B}\geq 0$, then $(A-\omega)_+-(B-\omega)_+\leq (A+\tilde{A}-\omega)_+ - (B+\tilde{B}-\omega)_+$. If $x\leq x_0< y$ or $x\leq y\leq x_0$, the implication is straightforward, following a similar split up.

The second statement follows from a similar reasoning: suppose there is a Poisson point at $(x_0,t)$ with weight $\omega$. If $x> x_0\geq 0$,
\begin{eqnarray*}
M^{\bar\nu}_t((x_0,x])-M^\nu_t((x_0,x]) & = & (M^{\bar\nu}_{t-}((x_0,x]) - \omega)_+ - (M^{\nu}_{t-}((x_0,x]) - \omega)_+\\
&\leq & M^{\bar\nu}_{t-}((x_0,x])-M^{\nu}_{t-}((x_0,x]).
\end{eqnarray*}
The inequality follows from the fact that $(A-c)_+-(B-c)_+\leq A_+-B_+$ whenever $c\geq 0$ and $A\geq B$. Since $M^{\bar\nu}_t([0,x_0])-M^\nu_t([0,x_0]) =  M^{\bar\nu}_{t-}([0,x_0])-M^{\nu}_{t-}([0,x_0])$, this shows that
\[ M^{\bar\nu}_t([0,x])-M^\nu_t([0,x])\leq M^{\bar\nu}_{t-}([0,x])-M^{\nu}_{t-}([0,x]).\]
Now suppose $x_0<0$ and $x\geq 0$. Note that under the condition on $\bar \nu$, we have that for all $s\geq 0$ and all $\eps>0$, $L_{\bar\nu}(-\eps,s)=L_{\nu}(-\eps,s)$, so
\[ M^{\bar\nu}_s([0,x])-M^\nu_s([0,x]) = L_{\bar\nu}(x,s)-L_\nu(x,s).\]
When $L_{\bar\nu}(x,t)$ does not use the weight at $(x_0,t)$, we know that $L_{\bar\nu}(x,t)=L_{\bar\nu}(x,t-)$ and that $L_{\nu}(x,t)\geq L_{\nu}(x,t-)$, which implies the desired result. If $L_{\bar\nu}(x,t)$ does use the weight at $(x_0,t)$, then it is not hard to see that $L_\nu(x,t)$ will also use the weight at $(x_0,t)$ (the longest path corresponding to $\bar \nu$ is always to the right of the path corresponding to $\nu$), which means that only the mass on the $x$-axis strictly to the left of $0$ is used, and therefore $M^{\bar\nu}_t([0,x])=M^\nu_t([0,x])$. Finally, when $x_0=x$ or $x_0>x$, we get that $M^{\bar\nu}_{t-}([0,x])=M^{\bar\nu}_t([0,x])$ and $M^{\nu}_{t-}([0,x])=M^\nu_t([0,x])$.
\hfill $\Box$\\

In other words, Proposition \ref{prop:coupling} tells us that the interacting fluid system is monotone: if one starts the fluid process with the same Poisson weights (basic coupling) and with ordered initial configurations, then the order is preserved for all $t\geq 0$. This coupled process is called the multi-class fluid system. The multi-class system is just a convention to describe a coupled process with ordered initial configurations (Ferrari \& Martin \cite{FeFo}).\\

\subsection{The multi-class invariant process} With Theorem \ref{thm:existray} in hands, for any countable $D\subseteq(\pi,3\pi/2)$, one can construct simultaneously a collection of equilibrium processes $\{\nu_\alpha:\alpha\in D\}$ by using the same Poisson weights on $\R\times\R_-$ and the Busemann functions $B_{\alpha}$. It turns out that this collection respects the order induced by the angles $\alpha\in D$. More precisely:
\begin{thm}\label{thm:multi-class}
If $\bar\alpha>\alpha$ then $\nu_{\bar\alpha}\geq\nu_{\alpha}$. In particular, for any countable subset  $\{\alpha_{i}\,:\,i\in\ZZ\}\subseteq D$, if one runs simultaneously (basic coupling) the interacting fluid processes on $\R\times\R_+$ with initial measures $(\nu_{\alpha_i}\,:\,i\in\ZZ)$ then, whenever $\alpha_i>\alpha_j$, $M^{\nu_{\alpha_i}}_t\geq M^{\nu_{\alpha_j}}_t$ for all $t\geq 0$, and
$$(M^{\nu_{\alpha_i}}_t\,:\,i\in \ZZ)\stackrel{\cal D}{=}(\nu_{\alpha_i}\,:\,i\in \ZZ)\,.$$
\end{thm}
\noindent{\bf Proof:}
Let $z'\geq z \in \R$. Let $\mm$ be the crossing point between $\varpi_{\bar\alpha}((z,0))$ and $\varpi_\alpha((z',0))$. Furthermore, denote $\cc$ as the coalescence point of the two $\alpha$-rays $\varpi_\alpha((z,0))$ and $\varpi_\alpha((z',0))$, and denote $\bar\cc$ as the coalescence point of the two $\bar\alpha$-rays $\varpi_{\bar\alpha}((z,0))$ and $\varpi_{\bar\alpha}((z',0))$. Then
\begin{eqnarray}
\nonumber \nu_{\bar\alpha}\big([z,z']\big)- \nu_{\alpha}\big([z,z']\big)&=& \Big\{L(\bar\cc,(z',0))-L(\bar\cc,(z,0))\Big\}\\
\nonumber&-&\Big\{ L(\cc,(z',0))-L(\cc,(z,0)) \Big\}\\
\nonumber&=& L(\bar\cc,(z',0))-\big\{L(\bar\cc,\mm)+L(\mm,(z',0))\big\}\\
\nonumber&+& L(\cc,(z,0))-\big\{L(\cc,\mm)+L(\mm,(z,0))\big\}\,\,\geq \,\,0\,.
\end{eqnarray}
Notice that the Busemann functions $B_\alpha$ are a function of the compound Poisson process $(\PP,\ww)$: $B_\alpha(\cdot,\cdot)=B_\alpha(\cdot,\cdot)(\PP,\ww)$ (here $\ww$ denotes the weights). Since $\pp+\PP$, the translated version of $\PP$, has the same distribution as $\PP$, we get that
\[ \{ B_{\alpha_i}(\cdot,\cdot)(\PP,\ww)\ :\ i\in \ZZ\} \stackrel{\cal D}{=} \{ B_{\alpha_i}(\cdot,\cdot)(\pp+\PP,\ww)\ :\ i\in \ZZ\}.\]
This shows time invariance for the Busemann multi-class process.

\hfill $\Box$\\

This result is also new in the classical Hammersley interacting system, where a different and explicit description of the invariant process with a finite number of classes is given in Ferrari \& Martin \cite{FeFo}.\\

\subsection{Law of large numbers for second-class particles} Proposition \ref{prop:coupling} can be used to define the notion of second-class particles. In the interacting fluid system we can define it analogously to the interacting particle case, with a slight adaptation due to the continuous weights. We start by changing $\nu$ into $\bar\nu$, by putting an extra weight $\eps>0$ in $0$, so
\[\bar\nu([0,x])=\nu([0,x])+\eps\,\,\,\mbox{ for }\,\,\,x\geq 0\,.\]
With this new process, and using the same Poisson weights, we define $M^{\bar\nu}_t$. Clearly,
\[ M^{\bar\nu}_t([0,x])\leq M^\nu_t([0,x]) + \eps.\]
Now define the location of the second class particle $X_\nu(t)$ as
\[ X_\nu(t) = \inf\{x\geq 0: M^{\bar\nu}_t(x)=M^\nu_t(x)+\eps\}\,.\]
By Proposition \ref{prop:coupling}, $X_\nu(t)$ is a non-decreasing function of $t$, meaning that the second class particle moves to the right. In fact, the extra mass $\eps$ will spread out, and our definition coincides with the rightmost point of this spread-out mass. This is a natural choice, since we will show that it does not depend on the total mass $\eps$, while for example the leftmost point does depend on $\eps$.

There is the following important connection between the longest path description and the second class particle. Let $\nu^+$ be the process defined by $\nu^+(x)=\nu(x)$ for $x\geq0$, and by $\nu^+(x)=-\infty$ for $x<0$. We also define the process $\nu^-$ by setting $\nu^-(x)=0$ for $x\geq0$, and $\nu^-(x)=\nu(x)$ for $x<0$. Then
\[ L_{\nu^+}(x,t) = \left\{ \begin{array}{ll}
\sup\{L((z,0),(x,t)) + \nu(z)\ :\ 0\leq z\leq x\} & \mbox{if}\ x\geq 0\\
-\infty & \mbox{if}\ x<0,
\end{array}\right.\]
and
\[ L_{\nu^-}(x,t) = \sup\{L((z,0),(x,t)) + \nu(z)\ :\ z<0\ \mbox{and}\ z\leq x\}.\]
Clearly,
$$L_\nu(x,t)=\max\big\{L_{\nu^+}(x,t),L_{\nu^-}(x,t)\big\}\,.$$
Now suppose $x\geq 0$. If $L_{\nu^+}(x,t)\geq L_{\nu^-}(x,t)$, there exists a longest path that does not use any weight of $\nu$ on $(-\infty,0)$. This means that if we add a weight $\eps>0$ in the origin, $L_{\bar\nu}(x,t)=L_{\bar\nu^+}(x,t)=L_\nu(x,t)+\eps$. Using Proposition \ref{prop:construction}, we see that this means that $M^{\bar\nu}_t(x)=M^\nu_t(x)+\eps$, so $X_\nu(t)\leq x$. If on the other hand we start with $X_\nu(t)\leq x$, we conclude that $M^{\bar\nu}_t(x)=M^\nu_t(x)+\eps$, using Proposition \ref{prop:coupling} and the fact that $M^\nu_t$ and $M^{\bar\nu}_t$ are right-continuous. This in turn means that $L_{\bar\nu}(x,t)=L_{\nu}(x,t)+\eps$, which is only possible if $L_{\nu^+}(x,t)\geq L_{\nu^-}(x,t)$. We have shown that
\begin{equation}\label{eq:secondclass}
\{ X_\nu(t)\leq x\} = \{L_{\nu^+}(x,t)\geq L_{\nu^-}(x,t)\}\,.
\end{equation}
Note that this can be rewritten as
\begin{equation}\label{eq:secondclassZ}
\{ X_\nu(t)\leq x\} = \{ Z_\nu(x,t) \geq 0\}.
\end{equation}
This means that the path of the second class particle corresponds to a competition interface, a fact well known for the totally asymmetric exclusion process (Ferrari \& Pimentel 2005). This allows us to show that the second class particle satisfies a strong law whenever $\nu^+$ and $\nu^-$ have asymptotic intensities. The proof of this does not use a coupling of two invariant versions of the fluid process, as is usual in the interacting particle case, but it uses the longest path description in a direct way. We would like to point out that in our general set-up, with random weights on the Poisson points, we do not have an equivalent of Burke's Theorem. This means that the time-reversed process is not a Hammersley interacting fluid system. Therefore, the path of a second class particle in general does not coincide in law with a longest path in the interacting fluid system, in contrast to the classical case, where the statement is true. However, we do have the following connection.
\begin{prop}\label{prop:second-class}
Assume that the distribution of $\nu$ is translation invariant. Then, for any $t\geq 0$, we have that
\[ X_{\nu}(t) - x \stackrel{\cal D}{=} -Z_\nu(x,t).\]
\end{prop}
\noindent{\bf Proof:} This follows almost immediately from \eqref{eq:secondclassZ}, since that equality can be rewritten as
\[ \{ X_{\nu}(t) - x  \leq h\} = \{ Z_\nu(x + h,t) \geq 0\}.\]
Now use translation invariance to see that
\[ Z_\nu(x + h,t) \stackrel{\cal D}{=} Z_\nu(x,t) + h.\]
Combining these two equations proves the proposition.

\hfill $\Box$\\

When we consider all $\alpha$-rays starting at the line $\R\times\{t\}$ and we move from left to right, $(X_{\nu_\alpha}(t),t)$ is the first point where the $\alpha$-ray passes the origin. It is tempting to think that the $\alpha$-ray starting at $(X_{\nu_\alpha}(t),t)$ actually passes through the origin, but this is false in general. In fact, after time $t$, most $\alpha$-rays will have coalesced with other rays, and the crossings with the $x$-axis will be quite far apart; we would conjecture they are order $t^{2/3}$ apart.

Proposition \ref{prop:second-class} allows us to use Theorem \ref{thm:shape} and Lemma \ref{lem:genexitcontrol} to prove a strong law for the second class particle in the case of $\nu_\alpha$. However, we are able to prove a strong law even for deterministic initial conditions that satisfy a density property:

\begin{thm}\label{thm:second-class}
Assume that
\begin{equation}\label{eq:second-int}
\lim_{x\to\infty}\frac{\nu\big(x)}{x}=\lim_{x\to-\infty}\frac{\nu(x)}{x}=\frac{\gamma}{2}\sqrt{\tan\alpha}\,.
\end{equation}
Then, with probability one,
$$\lim_{t\to\infty}\frac{X_\nu(t)}{t}=\frac{1}{\tan\alpha}\,.$$
\end{thm}
\noindent{\bf Proof:} It is enough to prove the statement for $\alpha=5\pi/4$. Suppose $\eps>0$ and $X_\nu(t)\leq t - 2\eps t$. Define $n=\lfloor t\rfloor$. Then for $t$ large enough, we have
\[X_\nu(n)\leq (1-2\eps)n + 1 + 2\eps\leq (1-\eps)n\,.\]
By  \eqref{eq:secondclassZ}, this implies that $Z_\nu((1-\eps)n,n)\geq 0$. In the same way as in Lemma \ref{eq:genexitcontrol}, this can happen only for finitely many $n\geq 1$, which gives
\[ \liminf_{t\to \infty} \frac{X_\nu(t)}{t} \geq 1\,.\]
Bounding the limit from above can be done using the analogous argument.

\hfill $\Box$\\

\section{Lattice last-passage percolation}\label{sec:TASEP}
In the lattice last-passage percolation model we have a collection $\{\omega_\zz\,:\,\zz\in\ZZ^2\}$ of i.i.d. non negative random variables indexed by lattice points $\zz\in\ZZ^2$. In this set up, one can define last-passage times for $\xx\leq\yy$ (we put an $\ell$ on the top to indicate that it refers to the lattice model) by maximizing over up-right paths connecting $\xx$ to $\yy$:
$$L^\ell(\xx,\yy):=\max_{\varpi\in\Pi(\xx,\yy)}\big\{\sum_{\xx'\in\varpi}\omega_{\xx'}\big\}\,.$$
Similary to the preivious case, one can start with a non-decreasing process $\nu^\ell=(\nu^\ell(z)\,:\,z\in\ZZ)$ and define the evolution by
$$M_t^{\nu^\ell}((x,y]):= L_{\nu^\ell}(y,t)-L_{\nu^\ell}(x,t)\,,\ \ \ \ \ (x,y\in\ZZ,\,\,t\geq 0)\,$$
where
$$L_{\nu^\ell}(x,t):=\sup_{z\leq x} \left\{ \nu^\ell(z) + L^\ell((z,0),(x,t))\right\}\,$$
($\nu^\ell$ must have a left density analogous to (\ref{eq:leftden})). For a survey in lattice last-passage percolation (and its connection with the totally asymmetric exclusion process) we address to \cite{Timo}.

If $\omega_\zz$ has an exponential distribution of parameter one, then the limit shape is given by $(\sqrt{x}+\sqrt{t})^2$ and a result similar to Theorem \ref{thm:existray} holds (Pimentel \& Ferrari \cite{FePi}). This allows us to construct Busemann functions $B^\ell$ for $\alpha\in(\pi,3\pi/2)$. The same method developed to prove Theorem \ref{thm:invariant} can be applied to this case. Since for exponential weights we also know the invariant measure, we have that:
\begin{thm}\label{thm:lattice}
For $\alpha\in(\pi,3\pi/2)$ let
$$\nu^\ell_\alpha(x)=B^\ell_\alpha\left((0,0),(x,0)\right)\mbox{ for }\,\,\, x\in\ZZ\,.$$
Then $\nu^\ell_\alpha$ is the unique ergodic process on $\ZZ$ that satisfies
$$M_t^{\nu^\ell_\alpha}\stackrel{\cal D}{=} \nu^\ell_\alpha\,\,\,\mbox{ for all }\,\,\,t\geq 0\,.$$
In particular, for any $x,y\in\ZZ$, $M_t^{\nu^\ell_\alpha}(x,y)$ is distributed like a sum of i.i.d. exponential random variables of intensity
\begin{equation}\label{eq:lattice}
\rho(\alpha):=\frac{\sqrt{\cos\alpha}}{\sqrt{\cos\alpha}+\sqrt{\sin\alpha}}\,.
\end{equation}
\end{thm}
We remark that $\rho(\alpha)$ is obtained by solving
$$(1+\sqrt{\tan\alpha})^2=\frac{1}{\rho}+\frac{\tan\alpha}{1-\rho}\,.$$
The left-hand side corresponds to the limit shape of the lattice model, in direction $(1,\tan\alpha)$, while the right-hand side corresponds to summing up, in the same direction, the expectation of the equilibrium process in the horizontal and vertical axis. (Compare this with \eqref{eq:gamma} and its proof.) The dual process on the vertical axis, denoted by $\nu_\alpha^{\ell*}$, is distributed like a sum of i.i.d. exponential weights of intensity $1-\rho(\alpha)$; see also Bal\'azs, Cator and Sepp\"al\"ainen \cite{BaCaSep}.

To develop our method in the model with general i.i.d. weights, one would need an inequality similar to (\ref{eq:curvature}) for the respective limit shape. This is, however, one of the most challenging problems in lattice last (and first) passage percolation models.

\section{Appendix}

\subsection{Proof of Proposition \ref{prop:Busemann}} The translation invariance of the underlying compound two-dimensional Poisson process,
$$\{\omega_\qq\,:\,\qq\in(\PP+\zz)\}\stackrel{\cal D}{=}\{\omega_\pp\,:\,\pp\in\PP\}\,,$$
implies \eqref{prop:Buse1}.

Anti-symmetry follows directly from the definition of the Busemann function. Now, by taking a coalescence point  $\cc$ between $\varpi_\alpha(\xx)$, $\varpi_\alpha(\yy)$ and $\varpi_\alpha(\zz)$ we have that
\begin{eqnarray}
\nonumber B_\alpha(\xx,\zz)&=&L(\cc,\zz)-L(\cc,\xx)\\
\nonumber&=&L(\cc,\zz)-L(\cc,\yy)+L(\cc,\yy)- L(\cc,\xx)\\
\nonumber&=&B_\alpha(\yy,\zz)+B_\alpha(\xx,\yy)\,,
\end{eqnarray}
which clearly shows additivity, and finishes the proof of \eqref{prop:Buse2}.

Define $S$ as the reflection in the diagonal $x=t$. Then $S(\PP)$ has the same distribution as $\PP$, and $\xx\leq \yy \Leftrightarrow S(x)\leq S(y)$. This shows that for all $\xx,\yy\in\R^2$,
\[ L(\xx,\yy) \stackrel{\cal D}{=} L(S(\xx),S(\yy)).\]
This clearly proves \eqref{prop:Busesym}.

Since $L(\cc,\yy)\geq L(\cc,\xx)$ whenever $\cc\leq\xx\leq\yy$ we have that $B_\alpha(\yy,\xx)\geq 0$ whenever $\xx\leq\yy$. Note that
\[ B_\alpha(\xx,\xx+(1,1))\geq L(\xx,\xx+(1,1)).\]
Clearly, $\E L(\xx,\xx+(1,1))>0$. Also,
\[\E B_\alpha(\xx,\xx+(1,1))=\E B_\alpha(\xx,\xx+(1,0)) + \E B_\alpha(\xx,\xx+(0,1)).\]
Now the first three properties prove that $\E B_\alpha(\xx,\yy)>0$.

To prove that its expected value is finite, without loss of generality, assume that $\yy=0$, that $\alpha=5\pi/4$, and that $\xx=(-1,-1)$. Let $\zz$ be the intersection point between $\varpi_{5\pi/4}(\0)$ and the one dimensional boundary of $\{\,\pp\in\R^2\,:\,\pp\leq\xx\,\}$. By taking a coalescence point $\cc$ for $\0$ and $\xx$ such that $\cc\leq \zz$, we have that
\begin{eqnarray}
\nonumber B_{5\pi/4}(\xx,\0)&=&L(\cc,\0)-L(\cc,\xx)\\
\nonumber &=& L(\cc,\zz)+L(\zz,\0)-L(\cc,\xx)\\
\nonumber &=& L(\zz,\0)- \{L(\cc,\xx)-L(\cc,\zz)\}\\
\nonumber&\leq&L(\zz,\0)\,.
\end{eqnarray}
Now assume that $\eta >0$ is a small constant and that $|\zz|\leq \eta r$. Then
$$L(\zz,\0)\leq \max\{L\left((-1,-\eta r),\0\right),L\left((-\eta r,-1),\0\right)\}\leq P(\eta r)\,$$
where $P(r)$ denotes the number of Poisson points in the strip composed of points $(x,t)$ such that $x\in[-1,0]$ and $t\in[-r,0]$, or $t\in[-1,0]$ and $x\in[-r,0]$. Hence
$$\P\left(L(\zz,\0)\geq r\right)\leq \P\left(|\zz|> \eta r\right)+\P\left(P(\eta r)\geq r\right)\,.$$
Denote by ${\rm Co}((-1,-1),\theta)$ the cone with axis the half-line starting at $\0$ through $(-1,-1)$, and angle $2\theta$. Lemma 2.3 of \cite{CPshape} ($\delta$-straightness of geodesics), states that for $\delta \in (0,1/4)$, $\theta=\pi/5$ and $M>0$ big enough, there exist $c_0,c_1>0$ such that for all $\pp\in {\rm Co}((-1,-1),\theta)$ and $|\pp|\geq M$,
\[\P\left(\{ \qq\in\R^2\ :\ \pp\in\varpi(\qq,\0)\}\not\subset {\rm Co}(\pp,c|\pp|^{-\delta})\right) \leq c_0e^{-c_1|\pp|^\kappa}\,.\]
If $|z|>\eta r$, then there must be $\pp\in \partial{\rm Co}((-1,-1),\theta)$ with $|\pp|>\eta r$, such that $\pp \in \varpi_{5\pi/4}(\0)$. However, for $r$ big enough, $\varpi_{5\pi/4}(\0)$ will leave the cone ${\rm Co}(\pp,c|\pp|^{-\delta})$, since it must follow the (negative) diagonal, implying
\[ \{ \qq\in\R^2\ :\ \pp\in\varpi(\qq,\0)\}\not\subset {\rm Co}(\pp,c|\pp|^{-\delta}).\]
Therefore, $\P(|\zz|>\eta r)$ is integrable in $r$ (for all $\eta$). On the other hand, by choosing $\eta$ small enough, one can make  $\P\left(P(\eta r)\geq r\right)$ integrable over $r>0$ as well. This finally implies that $\E L(\zz,\0)<\infty$, and consequently, $\E B_{5\pi/4}(\xx,\0)<\infty$.

By additivity and anti-symmetry, to prove cadlag, we can restrict our attention to $\xx,\yy,\pp=\0$, $\lambda\downarrow 0$ and $\qq=(1,0)$ so that we are varying in the horizontal direction close to the origin. For $\qq=(0,1)$ the argument is similar. Fix $\epsilon_0>0$ and assume that $B_\alpha(\0,(\epsilon_0,0))>\delta$. Then one can find $\cc_0$ such that the path $\varpi(\cc_0,(\epsilon_0,0))$ is disjoint from the path $\varpi(\cc_0,(0,0))$ and contains at least one Poisson point $\pp^0=(p^0_1,p^0_2) $ with $p^0_1\in(0,\epsilon_0)$. Define the set $A$ as the region enclosed by $\varpi(\cc_0,\0)$, the $x$-axis and the $y$-axis. Now take $\epsilon_1=p_1^0/2$ and assume that $B_\alpha(\0,(\epsilon_1,0))>\delta$. Then, similarly, one can find $\cc_1\geq \cc_0$ such that the path $\varpi(\cc_1,(\epsilon_1,0))$ is disjoint from the path $\varpi(\cc_1,(0,0))$ and contains at least one Poisson point $\pp^1=(p^1_1,p^1_2) $ with $p^1_1\in(0,\epsilon_1)$. But, now, we must have that $p^1\in A$. In this way, if $\liminf_{\epsilon\to 0}B_\alpha(\0,(\epsilon,0))>\delta$, then one could find infinitely many Poisson points within the bounded set $A$, which contradicts the fact that, a.s., the Poisson process is locally finite. \hfill $\Box$\\

\section{Acknowledgements}
Part of this work was done during our stay at the Institute Henri Poincare - Centre Emile Borel, attending the program Interacting Particle Systems, Statistical Mechanics and Probability Theory (September 5th - December 19th, 2008). Both authors wish to thank the organizers and the Institute for their hospitality and support during our stay there. Furthermore, we are grateful for careful reading and the many helpful suggestions by anonymous referees.

\end{document}